\documentclass{article}
\usepackage{amssymb, amsmath}
\usepackage{graphics}
\usepackage[dvips]{graphicx}
\usepackage{eepic}
\usepackage{epic}
\newtheorem{thm}{Theorem}[section]

\newtheorem{lmm}[thm]{Lemma}

\newcommand {\E} {\mathbb{E}}
\DeclareMathOperator {\var}{Var_{[0,1]}}

\newcommand {\p} {\mathbb{P}}
\newcommand {\Z} {\mathbb{Z}}
\newcommand {\N} {\mathbb{N}}
\newcommand {\R} {\mathbb{R}}

\newcommand {\ve} {\varepsilon}

\allowdisplaybreaks
\title{On permutations of lacunary series}
\author{Christoph \textsc{Aistleitner},\,
           Istv\'an \textsc{Berkes}
\,\, and Robert \textsc{Tichy}}
\begin{document}
%
\date{}
\maketitle

\begin{abstract}      
It is a well known fact that for periodic measurable $f$ and rapidly increasing $(n_k)_{k \geq 1}$  the sequence $(f(n_kx))_{k\ge 1}$
behaves like a sequence of independent, identically distributed random variables. For example, if
$f$ is a periodic Lipschitz function, then $(f(2^kx))_{k\ge 1}$ satisfies the central limit theorem, the law of
the iterated logarithm and several further limit theorems for i.i.d.\ random variables. Since an
i.i.d.\ sequence remains i.i.d.\ after any permutation of its terms, it is natural  to expect that the
asymptotic properties of lacunary series are also permutation-invariant. Recently, however,
Fukuyama (2009) showed that a rearrangement of the sequence $(f(2^kx))_{k\ge 1}$ can change
substantially its asymptotic behavior, a very surprising result. The purpose of the present paper
is to investigate this interesting phenomenon in detail and to give necessary and sufficient
criteria for the permutation-invariance of the CLT and LIL for $f(n_kx)$.

\vspace{4.2cm}\footnotesize{\noindent Authors:\\
Christoph Aistleitner, Graz University of Technology, Department for Analysis and Computational
Number Theory, Steyrergasse 30, 8010 Graz, Austria. Research supported by FWF grant S9603-N23. \newline e-mail:
\texttt{aistleitner@math.tugraz.at}\\
Istv\'an Berkes, Graz University of Technology, Institute of
Statistics, M\"unzgrabenstra{\ss}e 11, 8010 Graz, Austria. Research
supported by the FWF Doctoral Program on Discrete Mathematics (FWF DK W1230-N13), FWF grant S9603-N23 and OTKA grants K 67961 and K 81928. \newline e-mail: \texttt{berkes@tugraz.at}\\
Robert Tichy, Graz University of Technology, Department for Analysis and Computational Number Theory,
Steyrer\-gasse 30, 8010 Graz, Austria. Research supported by the FWF Doctoral Program on Discrete Mathematics (FWF DK W1230-N13) and FWF grant S9603-N23. \newline e-mail: \texttt{tichy@tugraz.at}
}
\end{abstract}

\newpage
\section{Introduction}

By classical results of Salem and Zygmund \cite{sz1947} and Erd\H{o}s and G\'al \cite{eg}, if $(n_k)_{k\ge 1}$ satisfies
the Hadamard gap condition
\begin{equation}\label{had}
n_{k+1}/n_k \ge q >1 \qquad (k=1, 2, \ldots)
\end{equation}
then $(\cos 2\pi n_kx)_{k\ge 1}$ obeys the central limit theorem and the law of the iterated logarithm, i.e.
\begin{equation}\label{cosclt}
\lim_{N\to\infty} \lambda \left\{x\in (0, 1):\sum_{k=1}^N
\cos 2\pi n_kx\le t\sqrt{N/2} \right\}=(2\pi)^{-1/2}  \int_{-\infty}^t
e^{-u^2/2}du,
\end{equation}
and
\begin{equation}\label{coslil}
\limsup_{N\to\infty} \frac{\sum_{k=1}^N
\cos 2\pi n_kx}{\sqrt{N\log\log N}}=1 \qquad \textup{a.e.}
\end{equation}
where $\lambda$ denotes the Lebesgue measure. Philipp and Stout \cite{phs} showed that under (\ref{had})
on the probability space $([0, 1], {\cal B}, \lambda)$ there exists a
Brownian motion process $\{W(t), ~t\ge 0\}$ such that
\begin{equation}\label{br}
\sum_{k=1}^{N} \cos 2\pi n_kx=W(N/2)+O\left(N^{1/2-\rho}\right)\qquad
\textup{a.e.}
\end{equation}
for some $\rho>0$. This relation implies not only the
CLT and LIL for $(\cos 2\pi n_kx)_{k \geq 1}$, but a whole class of  further limit theorems for  independent,
identically distributed random variables; for examples and discussion we refer to \cite{phs}.
Similar results hold for lacunary sequences $f(n_kx)$, where $f$ is a measurable function on $\mathbb R$ satisfying
the conditions
\begin{equation}\label{fcond}
f(x+1)=f(x), \qquad \int_0^1 f(x)\, dx=0, \qquad  \int_0^1 f^2(x)\, dx<\infty.
\end{equation}
For example, Kac \cite{ka} proved that if (\ref{fcond}) holds and $f$ is either a Lipschitz function or is of bounded variation in $(0, 1)$, then
$f(2^kx)$ satisfies the central limit theorem, more precisely
\begin{equation*}\label{fcltkac}
\lim_{N\to\infty} \lambda \left\{x\in (0, 1):\sum_{k=1}^N
f(2^kx)\le \sigma t\sqrt{N}\right\}=(2\pi)^{-1/2}  \int_{-\infty}^t
e^{-u^2/2}du,
\end{equation*}
where
$$
\sigma^2=\int_0^1 f^2(x)\, dx+2\sum_{k=1}^\infty \int_0^1 f(x) f(2^kx)~ dx.
$$
The corresponding LIL
\begin{equation*}\label{flil}
\limsup_{N\to\infty} \frac{\sum_{k=1}^N f(2^kx)}{\sqrt{2\sigma^2 N\log\log N}}=1 \qquad \textup{a.e.}
\end{equation*}
was proved by Maruyama \cite {ma} and Izumi \cite{iz} and the analogue of (\ref{br}) for $f(2^kx)$ was proved by Berkes and Philipp \cite{beph2}. These results show that lacunary subsequences of the trigonometric system and of the system $(f(nx))_{n\ge 1}$ behave like sequences of i.i.d.\ random variables and since an i.i.d.\ sequence remains i.i.d.\ after any permutation of its terms, it is natural to expect that limit theorems for $f(n_kx)$ are also permutation-invariant. Recently, however, Fukuyama \cite{ft2} showed the surprising fact that rearrangement of the sequence $\left(\{2^kx\}\right)_{k \geq 1}$, where $\{ \cdot \}$ denotes fractional part, can change the LIL behavior of its discrepancy. The purpose of this paper is to show that the permutation- invariance of $f(n_kx)$ is intimately connected with the number theoretic properties of $(n_k)_{k \geq 1}$ and we will give necessary and sufficient criteria for the rearrangement-invariant CLT and LIL for $f(n_kx)$ in terms of the number of solutions of the 
  Diophantine equation
$$
a_1 n_{k_1}+\ldots +a_p n_{k_p}=b, \qquad 1\le k_1,, \ldots, k_p \le N.
$$
In Section 2 we will formulate our results for the trigonometric system and in Section 3 for the system $f(nx)$. The proofs will be given in Section 4.

\section{The trigonometric system}

In this section we deal with permutation-invariance of lacunary trigonometric series.

\begin{thm}\label{trig} Let $(n_k)_{k \ge 1}$ be a sequence of positive integers  satisfying the Hadamard gap condition (\ref{had}) and let $\sigma: {\mathbb N}\to{\mathbb N}$ be a permutation of the positive integers. Then we have  \begin{equation}\label{permclt}
N^{-1/2}\sum_{k=1}^N \cos 2\pi n_{\sigma (k)}x \overset{\mathcal D}{\longrightarrow} {\mathcal N} (0, 1/2)
\end{equation}
and
\begin{equation}\label{permlil}
\limsup_{N\to\infty} \frac{\sum_{k=1}^N \cos 2\pi n_{\sigma (k)}x}{\sqrt{N\log\log N}}=1 \qquad \textup{a.e.}
\end{equation}
\end{thm}

Here we used probabilistic terminology and $\stackrel{\mathcal D}
\longrightarrow$  denotes convergence in distribution in the probability space $([0, 1], {\mathcal B}, \lambda)$.

Theorem \ref{trig} shows that for Hadamard lacunary trigonometric series the CLT and LIL are permutation-invariant without any additional number theoretic
condition on $n_k$. We note also that the unpermuted CLT (\ref{cosclt}) and LIL (\ref{coslil}) hold for $\cos 2\pi n_kx$ under the weaker gap condition
\begin{equation*}\label{ergap}
n_{k+1}/n_k\ge 1+ ck^{-\alpha},  \qquad 0\le \alpha<1/2
\end{equation*}
(see Erd\H{o}s \cite{er62} and Takahashi \cite{tak72,tak75}). However, as the next theorem shows (which is contained in \cite{aibeti10}), for having the permuted CLT and LIL the Hadamard gap condition (\ref{had}) in Theorem \ref{trig} cannot be weakened.

\begin{thm} For any positive sequence
$(\varepsilon_k)_{k \geq 1}$ tending to 0, there exists a sequence
$(n_k)_{k \geq 1}$ of positive integers satisfying
\begin{equation*}\label{epsilongap}
n_{k+1}/n_k\ge 1+\varepsilon_k, \qquad k\ge k_0
\end{equation*}
and a permutation $\sigma: {\mathbb N}\to {\mathbb N}$ of the
positive integers such that
\begin{equation*}\label{ell1}
a_N^{-1} \sum_{k=1}^N \cos 2\pi n_{\sigma(k)} x -b_N
\stackrel{\mathcal D}
\longrightarrow G\\
\end{equation*}
where $G$ is the Cauchy distribution with density $\frac{1}{\pi (1+x^2)}$,  $a_N=\sqrt{N/\log N}$ and $b_N=O(1)$.
Moreover, there exists a
permutation $\sigma: {\mathbb N}\to {\mathbb N}$ of the positive
integers such that
\begin{equation} \label{lilcount}
\limsup_{N \to \infty} \frac{\sum_{k=1}^N \cos 2 \pi n_{\sigma(k)}
x} {\sqrt{N  \log N}} >0 \quad \textup{a.e.}
\end{equation}
\end{thm}

For subexponentially growing $(n_k)_{k \ge 1}$ the permuted CLT and LIL are much harder problems and we do not have a precise charaterization. Our next theorem gives sufficient Diophantine conditions in this case.
Let us say that a sequence $(n_k)_{k \geq 1}$ of positive integers satisfies

\medskip
{\bf Condition ${\mathbf R}$}, if for any $p\ge 2$  and any fixed nonzero integers $a_1, \ldots a_p$ satisfying $|a_1|\le p, \ldots, |a_p|\le p$, the Diophantine equation
\begin{equation}\label{R}
a_1 n_{k_1} + \ldots +a_p n_{k_p}=0, \qquad k_1<\ldots <k_p
\end{equation}
has only a finite number of solutions.
If the number of solutions of (\ref{R}) is at most  $C\exp (p^\alpha)$ for some constants $C, \alpha>0$, we say that $(n_k)_{k \ge 1}$ satisfies condition ${\mathbf R^*}$. 

\begin{thm}\label{permcltlil}
If $(n_k)_{k \geq 1}$ is a sequence of positive integers satisfying
condition $\bf R$, then for any permutation $\sigma: {\mathbb N}\to{\mathbb N}$ of the positive integers we have the central limit theorem (\ref{permclt}). If $(n_k)_{k\ge 1}$ satisfies condition $\mathbf R^*$,
we also have the permuted LIL (\ref{permlil}).
\end{thm}

It is not easy to decide if an explicitly given subexponential sequence satisfies condition $\bf R$ or $\bf R^*$. A simple example is the sequence $(n_k)_{k\ge 1}$ consisting of the numbers of the form $p_1^{k_1} \cdots p_r^{k_r}$ \hbox{$(k_1, \ldots, k_r\ge 0)$}, where $\{p_1, \ldots, p_r\}$ is a finite set of coprime integers. In this case the validity of $\bf R^*$ follows from a recent deep bound, due to Amoroso and Viada \cite{am}, in the subspace theorem of W.\ Schmidt \cite{sch}. On the other hand, it is easy to obtain sequences satisfying condition $\bf R^*$ via random constructions. Let
$\omega_k = (\log k)^\alpha$, $\alpha>0$ and let $n_k$ be independent random variables, defined on some probability space, such that $n_k$ is uniformly distributed over the integers of the interval  $[a(k-1)^{\omega_{k}}, ak^{\omega_k}]$, where $a$ is a large constant. Then with probability one, $(n_k)_{k\ge 1}$ satisfies condition $\bf R^*$. 
The so constructed sequence $(n_k)_{k\ge 1}$ grows much slower than exponentially, in fact its growth speed is only slightly faster than polynomial. In contrast to the CLT (\ref{cosclt}) which can hold even for sequences $(n_k)$ with $n_{k+1}-n_k=O(1)$ (see Fukuyama \cite{fu2010}) we do not know if there are polynomially growing sequences $(n_k)_{k\ge 1}$ satisfying the permutation-invariant CLT and LIL.

\section{The system $f(k x)$}\label{fsection}

Let $f$ be a measurable function satisfying
\begin{equation*} \label{fcond*}
f(x+1)=f(x), \qquad \int_0^1 f(x)~dx=0, \qquad \int_0^1 f^2(x)~dx < \infty
\end{equation*}
and let $(n_k)_{k \geq 1}$ be a sequence of integers satisfying the Hadamard
gap condition (\ref{had}). The central limit theorem for $f(n_kx)$
has a long history. Kac \cite{ka} proved that if $f$ is Lipschitz
continuous then $f(n_kx)$ satisfies the CLT for $n_k=2^k$ and not
much later Erd\H{o}s and Fortet (see \cite[p.\ 655]{kac1949})
showed that this is not any more valid if $n_k=2^k-1$. Gaposhkin
\cite{gap1966} proved that $f(n_kx)$ obeys the CLT if
$n_{k+1}/n_k\to\alpha$ where $\alpha^r$ is irrational for $r=1, 2
\ldots$ and the same holds if all the fractions $n_{k+1}/n_k$ are
integers. To formulate more general criteria, let us say that a
sequence $(n_k)_{k \geq 1}$ of positive integers satisfies\\

{\bf Condition ${\mathbf D}_2$}, if for any fixed  nonzero
integers $a, b, c$ the number of solutions of the Diophantine
equation
\begin{equation}\label{gap}
an_k+bn_l=c, \qquad k, l\ge 1
\end{equation}
is bounded by a constant $K(a, b)$, independent of $c$.\\

{\bf Condition ${\mathbf D}^{*}_2$} (strong $\mathbf{D}_2$), if
for any fixed integers $a\ne 0$, $b\ne 0$, $c$ the number of
solutions of the Diophantine equation (\ref{gap}) is bounded by a
constant $K(a, b)$, independent of $c$, where for $a=b,c=0$ we require
also $k\ne l$.\\

Condition ${\mathbf D}_2$ is a variant of Sidon's $\bf B_2$ condition,
where the bound for the number of solutions of (\ref{gap}) is assumed only for $a=b=1$.  Gaposhkin \cite{gap1970} proved that under minor smoothness assumptions on $f$, condition ${\mathbf D}_2$ implies the CLT for $f(n_kx)$. Recently,  Aistleitner and Berkes \cite{aibe} proved
that the CLT holds for $f(n_kx)$ provided for any fixed nonzero integers $a,b,c$ the number of solutions $(k,l)$ of
$$
a n_k + b n_l = c,\qquad 1 \leq k,l \leq N
$$
is $o(N)$, uniformly in $c$, and this condition is also necessary. This criterion settles the CLT problem for $f(n_kx)$, but, as we noted, the validity of the CLT does not imply permutation-invariant behavior of $f(n_kx)$. The purpose of this section is to give a precise characterization for the CLT and LIL for permuted sums $\sum_{k=1}^N
f(n_{\sigma(k)}x)$.\\

 Our first result gives a necessary and sufficient condition
for permuted partial sums $\sum_{k=1}^N f(n_{\sigma(k)}x)$ to have
only Gaussian limit distributions and gives precise criteria this
to happen for a specific permutation $\sigma$.

\begin{thm}\label{th3}
Let $(n_k)_{k \geq 1}$ be a sequence of positive integers satisfying the
Hadamard gap condition (\ref{had}) and condition $\mathbf{D}_2$.
Let $f$ satisfy (\ref{fcond}) and let $\sigma$ be a permutation of
$\mathbb N$.
Then $ N^{-1/2}\sum_{k=1}^N f(n_{\sigma (k)}x)$ has a limit
distribution iff
\begin{equation}\label{gamma}
\gamma= \lim_{N\to\infty} N^{-1}\int_0^1 \left(\sum_{k=1}^N
f(n_{\sigma (k)}x)\right)^2 dx \
\end{equation}
exists, and then
\begin{equation}\label{fclt}
 N^{-1/2}\sum_{k=1}^N f(n_{\sigma (k)}x) \stackrel{\mathcal D}
\longrightarrow N(0, \gamma).
\end{equation}
If $\gamma=0$, the limit distribution is degenerate.
\end{thm}

Theorem \ref{th3} is best possible in the following sense:

\begin{thm} \label{th3bp}
If condition $\mathbf{D}_2$ fails, there exists a permutation
$\sigma$ such that the normed partial sums in (\ref{fclt}) have a
nongaussian limit distribution.
\end{thm}

In other words, under the Hadamard gap condition and condition
$\mathbf{D}_2$,  the limit distribution of $N^{-1/2}\sum_{k=1}^N
f(n_{\sigma (k)} x)$ can only be Gaussian, but the variance of the
limit distribution depends on the constant $\gamma$ in
(\ref{gamma}) which, as simple examples show, is not
permutation-invariant. For example, if $n_k=2^k$ and $\sigma$ is
the identity  permutation, then (\ref{gamma}) holds with
\begin{equation*}\label{kac}
\gamma=\gamma_f=\int_0^1 f^2(x)dx +2\sum_{k=1}^\infty \int_0^1
f(x)f(2^kx)dx
\end{equation*}
(see Kac \cite{ka}). Using an idea of Fukuyama \cite{ft2}, one can
construct permutations $\sigma$ of $\mathbb N$ such that
\begin{equation*} \label{varperm}
\lim_{N\to\infty} \frac{1}{N} \int_0^1 \left(\sum_{k=1}^N
f(n_{\sigma(k)} x)\right)^2 dx =\gamma_{\sigma, f}
\end{equation*}
with $\gamma_{\sigma, f} \ne \gamma_f$. Actually, the set of possible values $\gamma_{\sigma, f}$ belonging to all
permutations $\sigma$ contains the interval $I_f=[\|f\|^2,\gamma_f]$ and it is equal to this interval provided the Fourier
coefficients of $f$ are nonnegative. For general $f$ this is
false (for details cf.\  \cite{hlp}).\\

Under the slightly stronger condition $\mathbf{D}_2^{*}$ we get

\begin{thm}\label{co1} Let $(n_k)_{k \geq 1}$ be a sequence of positive integers satisfying the Hadamard gap condition (\ref{had}) and condition $\mathbf{D}_2^{*}$. Let $f$ satisfy (\ref{fcond}) and let $\sigma$ be a permutation of $\mathbb N$. Then the central limit theorem (\ref{fclt}) holds with $\gamma=\|f\|^2$.
\end{thm}

We now pass to the problem of the LIL.

\begin{thm} \label{thlil}
Let $(n_k)_{k \geq 1}$ be a sequence of positive integers satisfying the
Hadamard gap condition (\ref{had}) and condition $\mathbf{D}_2$.
Let $f$ be a measurable function satisfying (\ref{fcond}), let
$\sigma$ be a permutation of $\mathbb N$ and assume that
\begin{equation*}\label{gammastrong}
\gamma = \lim_{N \to \infty} N^{-1} \int_0^1 \left(\sum_{k=1}^N f(n_{\sigma (k)}x)\right)^2 dx
\end{equation*}
for some $\gamma\ge 0$. Then we have
\begin{equation}\label{lilperm3}
\limsup_{N\to \infty} \frac{\sum_{k=1}^N
f(n_{\sigma (k)} x)}{\sqrt{2N\log\log N}} =\gamma ^{1/2}\quad\textup{a.e.}
\end{equation}
If instead of condition $\mathbf{D}_2$ we assume $\mathbf{D}_2^{*}$, in (\ref{lilperm3}) we have $\gamma=\|f\|$.

\end{thm}

Theorem \ref{th3bp} shows that condition $\bf D_2$ in Theorem \ref{th3} is best possible: if there exist nonzero integers $a,b$ and a sequence $(c_m)_{m \geq 1}$ of integers such that the number of solutions $(k,l),~k \neq l$ of
$$
a n_k + b n_l = c_m
$$
tends to infinity as $m \to \infty$, the CLT for $f(n_{\sigma(k)} x)$ fails for a suitable permutation $\sigma$ and a trigonometric polynomial $f$. We cannot prove the analogous statement in Theorem \ref{thlil}. However, the conclusion of Theorem \ref{thlil} fails to hold for appropriate $f$ and $\sigma$ provided there exist nonzero integers $a,b$, and a fixed integer $c$, such that the Diophantine equation (\ref{gap}) has infinitely many solutions $(k,l), k \neq l$.

The previous results describe quite precisely the permutation-invariant CLT and LIL under the Hadamard gap condition (\ref{had}). If $(n_k)_{k \ge 1}$ grows faster than exponentially, i.e.
$$
n_{k+1}/n_k\to\infty
$$
then condition  $\mathbf{D}_2^*$ is satisfied, and thus so are the permutational-invariant CLT and LIL. In \cite{aibeti11} we have shown that given a 1-periodic Lipschitz $\alpha$ function $(0<\alpha\le 1)$ satisfying $\int_0^1 f(x)\, dx=0$, under the slightly stronger gap condition
$$
\sum_{k=1}^\infty (n_k/n_{k+1})^\alpha<\infty, \qquad 0<\alpha<1
$$
there exists a sequence $(g_k (x))_{k \geq 1}$ of measurable functions on $(0, 1)$, i.i.d.\ in the probabilistic sense, such that
$$
\sum_{k=1}^\infty |f(n_kx)-g_k(x)|<\infty \qquad \textup{a.e.}
$$
This implies not only the CLT and LIL for any rearrangement of $(f(n_kx))_{k\ge 1}$, but also the permutation-invariance of practically all asymptotic properties of $(f(n_kx))_{k \geq 1}$.

\section{Proofs}

To keep our paper at reasonable length, we will give here the proofs only for the CLT case. The arguments for the LIL use similar ideas, but they are considerably more complicated and will be given in a subsequent paper. As the proof of Theorem \ref{th3} will show, in the case when $f$ is a trigonometric polynomial of degree $d$, it suffices to assume condition $\mathbf{D}_2$ with coefficients $a,b$ satisfying  $|a|\leq d, |b| \leq d$.
In particular, in the trigonometric case $f(x)=\cos 2\pi x$ it suffices  to allow only coefficients $\pm 1$
in condition $\mathbf{D}_2$, i.e.\ condition $\mathbf{D_2}$ reduces to Sidon's $\mathbf{B_2}$ condition mentioned
in Section \ref{fsection}. This condition is satisfied for any Hadamard lacunary sequence $(n_k)_{k\ge 1}$
(see e.g.\ Zygmund \cite[pp.~203-204]{zt}) and thus Theorem \ref{trig} is contained in the proof of Theorem \ref{th3}.\\

Theorem \ref{co1} follows from Theorem \ref{th3} and the following
\begin{lmm} \label{lemmavarf2}
Let $(n_k)_{k \geq 1}$ be a lacunary sequence of positive integers satisfying condition $\mathbf{D}_2^{*}$, and let $f$ be a function satisfying (\ref{fcond}). Then for any permutation $\sigma$ of $\N$
$$
\lim_{N \to \infty} N^{-1} \int_0^1 \left( \sum_{k=1}^N f(n_{\sigma(k)} x)\right)^2~dx = \|f\|^2.
$$
\end{lmm}
\emph{Proof of Lemma \ref{lemmavarf2}}: For the simplicity of writing we assume that $f$ is an even function, i.e. the Fourier series of $f$ is of the form
$$
f \sim \sum_{j=1}^\infty a_j \cos 2 \pi j x.
$$
Let $\ve >0$ be arbitrary. We choose $d > \ve^{-1}$ and define
$$
p(x) = \sum_{j=1}^d a_j \cos 2 \pi j x, \qquad r(x)=\sum_{j=d+1}^\infty a_j \cos 2 \pi j x.
$$
Since by assumption $\var f < \infty$ we have
\begin{equation*} \label{sizefkoeff}
|a_j| = \mathcal{O} \left(j^{-1}\right), \qquad j \to \infty
\end{equation*}
(cf. Zygmund \cite[p.~48]{zt}), and therefore
\begin{equation} \label{rnorm}
\|r\| = \sum_{j=d+1}^\infty a_j^2/2 \ll \sum_{j=d+1}^\infty j^{-2} \ll d^{-1} \ll \ve
\end{equation}
By a classical norm inequality for lacunary series we have
\begin{equation} \label{normin}
\left\|\sum_{k=1}^N r(n_{\sigma(k)} x) \right\| \ll d^{-1/2} \sqrt{N},
\end{equation}
where the implied constant depends only on $f$ and the growth factor $q$ (see e.g.\ \cite[Lemma 2.1]{aibe}). Thus we have
$$
\left\| \sum_{k=1}^N r(n_{\sigma(k)} x) \right\| \ll \sqrt{\ve N}.
$$
By the orthogonality of the trigonometric system and Minkowski's inequality we have
\begin{eqnarray*}
& & \left\| \sum_{k=1}^N f(n_{\sigma(k)} x) \right\| \leq \left\| \sum_{k=1}^N p(n_{\sigma(k)} x) \right\| + \left\| \sum_{k=1}^N r(n_{\sigma(k)} x) \right\| \\
& \leq & \left(  \sum_{j=1}^d \sum_{k=1}^N \frac{a_j^2}{2} + \underbrace{\sum_{k_1,k_2=1}^N \sum_{j_1,j_2 = 1}^d}_{(j_1,k_1) \neq (j_2,k_2)} \frac{a_{j_1}a_{j_2}}{2} ~\mathbf{1}(j_1 n_{\sigma(k_1)} = j_2 n_{\sigma(k_2)}) \right)^{1/2} + \mathcal{O} \left( \sqrt{\ve N} \right).
\end{eqnarray*}
Since by condition $\mathbf{D}_2^{*}$ the number of nontrivial solutions of $j_1 n_{k_1} - j_2 n_{k_2}=0$ is bounded by a constant (where we can choose the same constant for all finitely many possible values of $j_1,j_2$), we get
\begin{equation} \label{ersn}
\left\| \sum_{k=1}^N f(n_{\sigma(k)} x) \right\| \leq \sqrt{N} \|p\| + \mathcal{O}(\sqrt{\ve N}).
\end{equation}
A lower bound can be found in exactly the same way, and since the implied constant in (\ref{ersn}) does not depend on $\ve$ (and $d$, resp.),
we obtain, utilizing (\ref{rnorm}),
$$
\left| \, \left\|\sum_{k=1}^N f(n_{\sigma(k)} x) \right\| - \sqrt{N} \|f\| \right| \ll \sqrt{\ve N},
$$
where the constant implied by $\ll$ is independent of $\ve$. Since $\ve$ can be chosen arbitrarily small, this proves the lemma. \qquad $\square$\\

The proof of Theorem \ref{permcltlil} is implicit in our paper \cite{hlp}. Let us say that an increasing sequence $(n_k)_{k\ge 1}$ of positive integers satisfies Condition $\mathbf R^{**}$ if for any $p\ge 2$ the number of nondegenerate solutions of the Diophantine equation
\begin{equation}\label{R**}
\pm  n_{k_1} \pm \ldots \pm  n_{k_p}=0, \qquad k_1, \ldots, k_p \ge 1
\end{equation}
is at most $C \exp (p^\alpha)$ for some constants $C, \alpha>0$ (here a solution of (\ref{R**}) is called nondegenerate if no subsum of the left hand side equals 0). As an analysis of the proof of the main theorem in \cite{hlp} shows, the proof uses only condition $\mathbf R^{**}$, and
by collecting the terms in (\ref{R**}) with equal indices shows that condition  $\mathbf R^{*}$ implies condition $\mathbf R^{**}$.

To prove the remark made at the end of Section 2, let
$\omega_k = (\log k)^\alpha$ and let $n_k$, $k=1, 2, \ldots$ be independent random variables such that $n_k$ is uniformly distributed over the integers of the interval $I_k = [a(k-1)^{\omega_{k-1}}, ak^{\omega_k}]$.  Note that the length of $I_k$ is at least $a\omega_k (k-1)^{\omega_k-1} \ge a\omega_1$ for $k=2, 3, \ldots$ and equals $a$ for $k=1$ and thus choosing $a$ large enough, each $I_k$ contains at least one integer. We claim that, with probability 1, the sequence $(n_k)_{k\ge 1}$ satisfies
condition $\mathbf{R^*}$.  To see this,
set $\eta_k=\frac{1}{2}\omega_k^{1/2}$, then
\begin{equation} \label{etak}
(2k)^{\eta_k^2+2 \eta_k} \leq (2k)^{\omega_k/2}
\le k^{-2}|I_k| \qquad \text{for} \ k\ge k_0
\end{equation}
since, as we noted, $|I_k|\ge a \omega_k (k-1)^{\omega_k-1}\ge
(k/2)^{\omega_k-1}$ for large $k$. Let $k\ge 1$ and
consider the numbers of the form
\begin{equation}\label{const}
(a_1 n_{i_1}+\ldots +a_ s n_{i_s})/d
\end{equation}
where $1 \leq s \leq \eta_k$,  $1\le i_1, \ldots, i_s\le k-1$,
$a_1, \ldots a_s, d$ are nonzero integers with
$|a_1|, \ldots, |a_s|, |d|\le \eta_k$. Since the number of values in (\ref{const}) is at most
$(2k)^{\eta_k^2+2 \eta_k}$, (\ref{etak}) shows that
the probability that $n_k$ equals any of these
numbers is at most $k^{-2}$. Thus by the Borel-Cantelli lemma,
with probability 1 for $k\ge k_1$, $n_k$ will be different from
all the numbers in (\ref{const}) and thus the equation
$$a_1 n_{i_1}+\ldots +a_s n_{i_{s}}+ a_{s+1} n_k=0$$
has no solution with $1 \leq s \leq \eta_k$, $1\le i_1< \ldots <i_s\le k-1$, $0< |a_1|, \ldots, |a_{s+1}|\le \eta_k$. By monotonicity, the equation
\begin{equation}\label{const2}
a_1n_{i_1}+\ldots + a_s n_{i_s}=0, \qquad i_1<\ldots < i_s
\end{equation}
has no solutions provided the number of terms is at most
$\eta_k$, $0< |a_1|, \ldots, |a_s|\le \eta_k$ and $i_s\ge k$. Thus the number  of solutions of (\ref{const2}) where $s\le \eta_k$ and $0< |a_1|, \ldots, |a_s|\le \eta_k$ is at most $k^{\eta_k}\ll \exp(\eta_k^\beta)$ for some $\beta>0$. Since the sequence $[\eta_k]$ takes all sufficiently large integers $p$, the sequence $(n_k)_{k\ge 1}$ satisfies condition $\mathbf{R^*}$.

\medskip
It remains now to prove Theorems \ref{th3} and \ref{th3bp}, which will be done in the next two sections of the paper.

\subsection{Proof of Theorem \ref{th3}.} \label{sectth3proof}

In this section we give the proof of Theorem \ref{th3}. The proof of the main lemma (Lemma \ref{lemmaclt}) uses ideas of R{\'e}v{\'e}sz \cite{rt}.
We will need the following simple smoothing inequality.
\begin{lmm} \label{lemmaclt2} \label{lemmacf1dim}~Let $P_1,P_2$ be probability measures on $\R$, and write $p_1,
p_2$ for the corresponding characteristic functions.  Let
$P_1^* = P_1 \star H$, $P_2^* = P_2 \star H$,
where $H$ is a normal distribution with mean zero and variance $\tau^2$. Then for all $y>0,T>0,$
\begin{eqnarray*}
\left|P_1^*([-y,y]) - P_2^*([-y,y])\right| & \leq & y \int_{s \in
[-T,T]} |p_1(s)-p_2(s)| ~ds \\
& & + 4 y  (\tau^{-1}~\exp(- T^2 \tau^2/2)).
\end{eqnarray*}
\end{lmm}

\emph{Proof:~} Letting $h(s)=\exp(-\tau^2 s^2/2)$,  the characteristic functions of  $P_1^*$ and $P_2^*$ are $p_1^*=p_1 h$ and $p_2^*=p_2 h$
and thus for the densities $\rho_1$ and $\rho_2$ of $P_1^*$ and $P_2^*$ we have, letting $w(s)=p_1(s)-p_2(s)$,
\begin{eqnarray*}
\left| \rho_1(u) - \rho_2(u)\right| & \leq & (2 \pi)^{-1} \left| \int_{\R} e^{-isu} \left(p_1^*(s)-p_2^*(s)\right)~ds\right| \\
& \leq & (2\pi)^{-1} \int_{s\in[-T,T]} |w(s)| ds + (2 \pi)^{-1} ~2 \int_{s \not\in [-T,T]} |h(s)|ds.
\end{eqnarray*}
Thus
\begin{eqnarray*}
\left|P_1^*([-y,y]) - P_2^*([-y,y])\right| & \leq & \int_{[-y,y]} |\rho_1(u)-\rho_2(u)|du \\
& \leq & y \int_{s \in [-T,T]} |w(s)| ~ds + y \int_{s \not\in[-T,T]}|h(s)|~ds,
\end{eqnarray*}
and the statement of the lemma follows from
\begin{eqnarray*}
\int_{s\not\in[-T,T]} |h(s)| ~ds = 2 \int_T^\infty e^{-\tau^2 s^2/2}~ds \leq & 4 \tau^{-1} e^{-\tau^2 T^2/2}. \qquad \square
\end{eqnarray*}

Let $f$ be a function satisfying (\ref{fcond}) and $(n_k)_{k\ge 1}$ a sequence satisfying the Hadamard gap condition (\ref{had}) and condition $\mathbf{D}_2$, and let
$$
p(x) = \sum_{j=1}^d (a_j \cos 2 \pi j x + b_j \sin 2 \pi j x)
$$
be the $d$-th partial sum of the Fourier series of $f$ for some $d \geq 1$. By condition $\mathbf{D}_2$ we can find a number $C_1$ (depending on $d$) such that for any $a,b$ satisfying $0 < |a|\leq d,~0 < |b|\leq d$
\begin{equation} \label{b2const}
\# \left\{k_1,k_2 \geq 1:~a n_{k_1} - b n_{k_2}=c \right\} \leq C_1,
\end{equation}
for all $c \in \Z \backslash \{0\}$. We set
$$
\gamma_N = N^{-1/2} \left\| \sum_{k=1}^N f(n_{\sigma(k)} x) \right\| .
$$

\begin{lmm}\label{lemmaclt}
There exists an $N_0$, depending only on $p,d$ and the  constants implied by condition $\mathbf{D}_2$ (but not on $\sigma$) such that for any $N \geq N_0$ there exists a set $A \subset \{1,\dots,N\}$ with $\#A \geq N - N (\log_q (2d))/(\log N)^{1/2}$ such that
$$
\left| \E \left(\exp \left(i s N^{-1/2}\sum_{k \in A} p(n_{\sigma(k)} x) \right)\right) - e^{-s^2 \delta_N^2/2} \right| \leq 2 N^{-1/8}
$$
for all $s \in [-(\log N)^{1/8}, (\log N)^{1/8}]$, where
$$
\delta_N = N^{-1} \left\| \sum_{k \in A} p(n_{\sigma(k)} x) \right\|.
$$
\end{lmm}

{\bf Proof.}
For the simplicity of writing we assume that
$$
p(x)=\sum_{j=1}^d a_j \cos 2 \pi j x
$$
is an even function; the proof in the general case is similar.
Let $(\nu_k)_{1 \leq k \leq N}$ denote the sequence  $(n_{\sigma(k)})_{1 \leq k \leq N}$ arranged in increasing order. Put
\begin{eqnarray*}
\Delta & = & \left\{k \in \{1, \dots, N\}:~k \mod \lceil (\log N)^{1/2} + \log_q (2d) \rceil \not\in \{0, \dots, \lceil \log_q (2d)\rceil \} \right\} \\
\Delta^{(h)} & = & \left\{ k \in \Delta:~\frac{k}{\lceil (\log N)^{1/2} + \log_q 2d \rceil} \in [h,h+1)\right\}, \qquad h \geq 0,\\
A & = & \{1 \leq l \leq N:~ n_{\sigma(l)} \in (\nu_k)_{k \in \Delta} \},\\
\alpha(s) & = & \prod_{h \geq 0} \left(1+is N^{-1/2} \sum_{k \in \Delta^{(h)}} p(\nu_k x)\right).
\end{eqnarray*}
and
$$
\eta_N=N^{-1/2}\sum_{k \in \Delta} p(\nu_k x).
$$
The sets $\Delta^{(h)}$ are constructed in such a way that for $k_1 \in \Delta^{(h_1)},~k_2 \in \Delta^{(h_2)}$, where $h_1 < h_2$, we have
\begin{equation} \label{rela}
\nu_{k_2}/\nu_{k_1} > q^{k_2-k_1} \geq q^{\lceil \log_q (2d) \rceil} \geq 2d.
\end{equation}
Using
\begin{equation} \label{win}
e^{ix} = (1+ix) e^{-x^2/2+w(x)}, \qquad |w(x)| \leq |x|^3
\end{equation}
we get
\begin{eqnarray}
e^{i s \eta_N} & = & \prod_{h \geq 0} \exp \left( is N^{-1/2} \sum_{k \in \Delta^{(h)}} p(\nu_k x)\right) \label{e1} \\
& = & \alpha(s)~\exp\left(-(2N)^{-1}\sum_{h \geq 0} s^2 \left(\sum_{k \in \Delta^{(h)}} p(\nu_k x) \right)^2 \right) \label{e2} \\
& & ~\times~\exp\left(\sum_{h \geq 0} w\left(isN^{-1/2} \sum_{k \in \Delta^{(h)}} p(\nu_k x) \right)\right). \label{e3}
\end{eqnarray}
We have
\begin{eqnarray}
& & \quad \sum_{h \geq 0} \left(\sum_{k \in \Delta^{(h)}} p(\nu_k x) \right)^2 \label{alf1} \\
& = & \sum_{h \geq 0} \sum_{k_1,k_2 \in \Delta^{(h)}} \sum_{j_1,j_2=1}^d \frac{a_{j_1} a_{j_2}}{2} \left( \cos (2\pi(j_1 \nu_{k_1}+j_2 \nu_{k_2})x)+\cos (2\pi(j_1 \nu_{k_1}-j_2 \nu_{k_2})x) \right) \nonumber \\
& = & \sum_{h \geq 0} \sum_{k_1,k_2 \in \Delta^{(h)}} \sum_{j_1,j_2=1}^d \frac{a_{j_1} a_{j_2}}{2} \left( \cos (2\pi(j_1 \nu_{k_1}+j_2 \nu_{k_2})x)\right) \nonumber\\
& & ~+ \sum_{h \geq 0} \underbrace{\sum_{k_1,k_2 \in \Delta^{(h)}} \sum_{j_1,j_2=1}^d}_{j_1 \nu_{k_1}-j_2 \nu_{k_2} = 0} \frac{a_{j_1}a_{j_2}}{2} \nonumber\\
& & ~+ \sum_{h \geq 0} \underbrace{\sum_{k_1,k_2 \in \Delta^{(h)}} \sum_{j_1,j_2=1}^d}_{j_1 \nu_{k_1}-j_2 \nu_{k_2} \neq 0} \frac{a_{j_1} a_{j_2}}{2} \left( \cos (2\pi(j_1 \nu_{k_1}-j_2 \nu_{k_2})x) \right) \nonumber\\
& &  \quad = \delta_N^2 N + \beta(x), \label{alf2}
\end{eqnarray}
where
\begin{eqnarray*}
\delta_N & = & N^{-1} \left(\sum_{h \geq 0} \underbrace{\sum_{k_1,k_2 \in \Delta^{(h)}} \sum_{j_1,j_2=1}^d}_{j_1 \nu_{k_1}-j_2 \nu_{k_2} = 0} \frac{a_{j_1}a_{j_2}}{2} \right)^{1/2}
\end{eqnarray*}
and
\begin{eqnarray*}
\beta(x) & = & \sum_{h \geq 0} \sum_{k_1,k_2 \in \Delta^{(h)}} \sum_{j_1,j_2=1}^d \frac{a_{j_1} a_{j_2}}{2} \left( \cos (2\pi(j_1 \nu_{k_1}+j_2 \nu_{k_2})x)\right) \\
& & ~+ \sum_{h \geq 0} \underbrace{\sum_{k_1,k_2 \in \Delta^{(h)}} \sum_{j_1,j_2=1}^d}_{j_1 \nu_{k_1}-j_2 \nu_{k_2} \neq 0} \frac{a_{j_1} a_{j_2}}{2} \left( \cos (2\pi(j_1 \nu_{k_1}-j_2 \nu_{k_2})x) \right).
\end{eqnarray*}
Note that
$$
\delta_N^2 = N^{-1} \int_0^1 \left(\sum_{h \geq 0} \sum_{k \in \Delta^{(h)}} p(\nu_k x)\right)^2 dx = N^{-1} \int_0^1 \left( \sum_{k \in \Delta} p(\nu_k x)\right)^2 dx,
$$
since by (\ref{rela}) the functions
$$
\sum_{k \in \Delta^{(h_1)}} p(\nu_k x) \qquad  \textrm{and} \qquad \sum_{k \in \Delta^{(h_2)}} p(\nu_k x)
$$
are orthogonal for $h_1 \neq h_2$. From (\ref{e1}), (\ref{e2}), (\ref{e3}), (\ref{alf1}) and (\ref{alf2}) we conclude
\begin{eqnarray*}
e^{i s \eta_N} = \alpha(s) ~\exp\left(-\frac{s^2}{2} ~\left(\delta_N^2+\frac{\beta}{N}\right) + \sum_{h \geq 0} w\left(isN^{-1/2} \sum_{k \in \Delta^{(h)}} p(\nu_k x) \right) \right)
\end{eqnarray*}
and, writing
$$
w^*(s,x)=\sum_{h \geq 0} w\left(is N^{-1/2}\sum_{k \in \Delta^{(h)}} p(\nu_k x)\right),
$$
(\ref{win}) implies
\begin{eqnarray}
|w^*(s,x)| & \leq & N \left|sN^{-1/2} \sum_{k \in \Delta^{(h)}} p(\nu_k x)\right|^3 \nonumber\\
& \leq & |s|^3 ~\|p\|_\infty^3 ~(\log N)^{3/2}~ N^{-1/2}. \label{w2}
\end{eqnarray}
We note that
\begin{eqnarray}
|\alpha(s)| & \leq & \prod_{h \geq 0} \left( 1 + N^{-1}s^2 \left( \sum_{k \in \Delta^{(h)}} p(\nu_k x) \right)^2\right)^{1/2} \label{als1}\\
& \leq & \exp \left( (2N)^{-1}\sum_{h \geq 0} s^2 \left( \sum_{k \in \Delta^{(h)}} p(\nu_k x) \right)^2\right) \nonumber\\
& \leq & \exp \left( \frac{s^2}{2} ~\left(\delta_N^2 + \frac{\beta(x)}{N} \right)\right) \label{als2}
\end{eqnarray}
and obtain, for $\varphi(s) = \E e^{i s \eta_N}$, using $\E \alpha(s)=1$ and (\ref{als1}), (\ref{als2}),
\begin{eqnarray}
& & \left| \varphi(s) - e^{-s^2 \delta_N^2/2} \right| \nonumber\\
& = & \left| \E \left( \alpha(s) ~\exp\left(-\frac{s^2}{2} ~ \left(\delta_N^2+\frac{\beta(x)}{N}\right) + w^*(s,x) \right) \right) - e^{-s^2 \delta_N^2/2} \right| \nonumber\\
& = & \left| \E \left( \alpha(s) \left(\exp\left(-\frac{s^2}{2} ~ \left(\delta_N^2+ \frac{\beta(x)}{N} \right) + w^*(s,x) \right) - e^{-s^2 \delta_N^2/2}\right) \right) \right| \nonumber\\
& \leq & \left| \E \left( |\alpha(s)|~\left| \exp\left(-\frac{s^2}{2} ~ \left(\delta_N^2+ \frac{\beta(x)}{N} \right) + w^*(s,x) \right) - e^{-s^2 \delta_N^2/2}\right| \right) \right| \nonumber\\
& \leq & \E \left|e^{w^*(s,x)} -1\right| + \E \left|\exp\left(\frac{s^2}{2}~\frac{\beta(x)}{N}\right)- 1\right|. \nonumber\label{distb1}
\end{eqnarray}
If $|s| \leq (\log N)^{1/8}$, then by (\ref{w2})
$$
\|p\|_\infty^3 ~(\log N)^{15/8}~ N^{-1/2} \leq N^{-1/4} \qquad \textrm{for} \qquad N \geq N_1
$$
with $N_1$ depending only on $p$, and hence
\begin{equation*} \label{distb2}
\left|e^{w^*(s,x)} -1\right| \leq 2 N^{-1/4}, \qquad \textrm{for}~N \geq N_1.\\
\end{equation*}
On the other hand, the function $\beta(x)$ is a sum of at most
$$
2 \sum_{h \geq 0} d^2 \left| \Delta^{(h)}\right|^2 \leq 2 d^2 N(\log N)^{1/2}
$$
trigonometric functions, and the coefficient of each of this summands is bounded by $\max_{1 \leq j \leq d} a_j^2/2 \leq C_2$, where $C_2$ depends only on $p$ and $d$. Since by assumption $(n_k)_{k \geq 1}$ satisfies condition $\mathbf{D}_2$, there can be at most $4 d^2 C_1$ summands giving the same frequency (the constant $C_1$ is defined in (\ref{b2const})). This means, writing $\beta(x)$ in the form
$$
\beta(x) = \sum_{j=1}^\infty c_j \cos 2 \pi j x
$$
we have
\begin{equation} \label{beta1}
|c_j| \leq 4 d^2 C_1 C_2, \quad j \geq 1, \qquad \textrm{and} \qquad \|\beta\|_\infty \leq \sum_{j=1}^\infty |c_j| \leq 2 C_2 d^2 N\sqrt{\log N}
\end{equation}
Therefore
$$
\|\beta\|^2 = \sum_{j=1}^\infty c_j^2/2 \leq \max_{j \geq 1} |c_j| \sum_{j=1}^{\infty} |c_j| \leq (4 d^2 C_1 C_2)^2  ~2 C_2 d^2 N \sqrt{\log N},
$$
and
$$
\p (|\beta|> N^{2/3} ) \leq N^{-1/4} \qquad \textrm{for} \qquad N \geq N_2,
$$
where $N_2$ depends only on $p,~d,~C_1$. Hence by (\ref{beta1})
\begin{eqnarray*}
& & \E \left|\exp\left(\frac{s^2}{2}~\frac{\beta(x)}{N}\right)- 1\right| \\
& \leq & \left|\exp\left(\frac{s^2}{2}~\frac{\|\beta\|_\infty}{N}\right)- 1\right| N^{-1/4} + \left|\exp\left(\frac{s^2}{2}~\frac{N^{2/3}}{N}\right)- 1\right|,
\end{eqnarray*}
and, assuming
$$
|s| \leq (\log N)^{1/8},
$$
we get
\begin{eqnarray}
& & \E \left|\exp\left(\frac{s^2}{2}~\frac{\beta(x)}{N}\right)- 1\right| \nonumber\\
& \leq & \left|\exp\left(C_2 d^2 (\log N)^{3/4} \right)- 1\right| N^{-1/4} + \left|\exp\left(\frac{(\log N)^{1/4} N^{2/3}}{N}\right)- 1\right| \nonumber\\
& \leq & N^{-1/8} \qquad \textrm{for} \qquad N \geq N_3, \nonumber\label{distb3}
\end{eqnarray}
where $N_3$ depends on $d,p,C_1$. Combining these estimates, we have
\begin{equation} \label{comb1}
|\varphi(s)-e^{s^2 \delta_N^2/2}| \leq 2 N^{-1/8} \quad \textrm{for} \quad N \geq N_4
\end{equation}
where $N_4$ also depends on $p,d,C_1$. This proves Lemma \ref{lemmaclt}. \quad $\square$\\

\noindent
\emph{Proof of Theorem \ref{th3}.~} Let $y$ and $\varepsilon>0$ be given. We choose $d > \varepsilon^{-2}$ and write $p$ and $r$, respectively, for the $d$-th partial sum and $d$-th remainder term of the Fourier series of $f$. Assume that $N \geq N_0$ with the $N_0$ in Lemma \ref{lemmaclt}, and let $A$ be the set in Lemma \ref{lemmaclt}. By Lemma \ref{lemmaclt} we have
$$
|\varphi_N(s)-e^{s^2 \delta_N^2/2}| \leq 2 N^{-1/8},
$$
where
\begin{eqnarray*}
\varphi_N(s) & = & \E e^{i s \eta_N}, \qquad \eta_N = \sum_{k \in \Delta} p(n_{\sigma(k)}x), \\
\delta_N^2 & = & N^{-1} \int_0^1 \left( \sum_{k \in A} p(n_{\sigma(k)} x)\right)^2 dx.
\end{eqnarray*}
We recall that
$$
\gamma_N^2 = N^{-1} \int_0^1 \left( \sum_{k=1}^N f(n_{\sigma(k)} x)\right)^2 dx.
$$
By (\ref{normin}) we have
\begin{equation} \label{fnorm}
\left\|\sum_{k=1}^N r(n_{\sigma(k)} x) \right\| \leq C_4 d^{-1/2} \sqrt{N},
\end{equation}
where $C_4$ depends on $f,q$. By Minkowski's inequality, and using (\ref{fnorm}) for our choice of $d>\varepsilon^{-2}$,
$$
\gamma_N \sqrt{N} \leq \left\|\sum_{k=1}^N p(n_{\sigma(k)} x) \right\| + \left\|\sum_{k=1}^N r(n_{\sigma(k)} x) \right\| \leq \left\|\sum_{k=1}^N p(n_{\sigma(k)} x) \right\| + C_4 \varepsilon \sqrt{N},
$$
and, since $\left\|\sum_{k=1}^N p(n_{\sigma(k)} x) \right\| \leq C_5 \sqrt{N}$ (for $C_5$ depending on $f,q$, again by \cite[Lemma 2.1]{aibe}), we obtain
$$
\gamma_N^2 \leq N^{-1} \left\|\sum_{k=1}^N p(n_{\sigma(k)} x) \right\|^2 + C_6 \varepsilon
$$
($C_6$ depends on $f,q$). Since $N - \#A \leq N(\log_q (2d))/\sqrt{\log N}$ we have
\begin{eqnarray}
\gamma_N^2 - \delta_N^2 & \leq & C_6 \varepsilon + 2 N^{-1} \int_0^1 \sum_{k_1,k_2 \in \{1,\dots,N\},~k_2 \not\in \Delta} p(\nu_{k_1} x) p(\nu_{k_2} x) ~dx \nonumber\\
& \leq & C_6 \varepsilon + 2 N^{-1} \|p\|^2 \sum_{k_2 \in \{1,\dots,N\}, k_2 \not\in \Delta} \# \{k_1 \in \{1, \dots, N\}:~k_1/k_2 \in [1/d,d]\} \nonumber\\
& \leq & C_6 \varepsilon + 2 N^{-1} \|p\|^2 N (\log N)^{-1/2} \lceil \log_q 2d \rceil^2 \nonumber\\
& \leq & C_7 \varepsilon \qquad \textrm{for} \qquad N \geq N_5, \label{n5est} 
\end{eqnarray}
where $C_7$ depends on $f,q$ and $N_5$ depends on $f,q,d$.  We assume $s\leq \varepsilon^{-1/2}$, and, in view of (\ref{comb1}), we arrive at
$$
|\varphi_N(s) - e^{s^2 \gamma_N^2/2}| \leq 2 N^{-1/8} + |e^{s^2 \gamma_N^2/2}-e^{s^2 \delta_N^2/2}| \leq C_8 \varepsilon
$$
for $N \geq N_6$, where $C_8$ depends on $f,q$ and $N_6$ depends on $f,q,d$. The proof of (\ref{n5est}) also shows
$$
\left\| \sum_{1 \leq k \leq N, ~k \not\in A} p(n_{\sigma(k)} x) \right\| \leq C_9 \varepsilon, \qquad N \geq N_7, 
$$
and thus
\begin{equation} \label{peste}
\p \left( \left| \sum_{1 \leq k \leq N, ~k \not\in A}  p(n_{\sigma(k)}x) \right| > \ve^{1/6} \right) \leq C_9 \ve^{2/3}, \qquad N \geq N_7,
\end{equation}
for $C_9$ depending on $f,q$ and $N_7$ depending on $f,q,d$.\\

We use Lemma \ref{lemmaclt2} with
$$
p_1(s)=\varphi(s), \qquad p_2(s) = e^{s^2 \gamma_N^2/2}, \qquad T=\varepsilon^{-1/2}, \qquad \tau=\varepsilon^{1/3}
$$
to get (using the notation from this lemma) for all $y>0$
$$
\left| P_1^* ([-y,y]) - P_2^*([-y,y]) \right| \leq y \int_{s \in [-T,T]} C_8 \ve ~ds + 4y \varepsilon^{-1/3} e^{-\ve^{-1/3}/2}
$$
for sufficiently large $N$ (depending on $f,q,d$), provided $\ve$ is sufficiently small, which  we can assume.
Thus, if $P_1,P_2$ are the measures corresponding to $p_1,p_2$, we get for sufficiently large $N$ (depending on $f,q,d$), using (\ref{fnorm}) and (\ref{peste}),
\begin{eqnarray}
& & \p \left( N^{-1/2} \sum_{k=1}^N f(n_{\sigma(k)}x) \in [-y,y]\right) \nonumber\\
& \leq & P_1 ([-y-2 \ve^{1/6},y+2 \ve^{1/6}]) + \p \left( \left| N^{-1/2} \sum_{k=1}^N r(n_{\sigma(k)}x) \right| > \ve^{1/6} \right)\nonumber\\ 
& & ~ + \p \left( \left| \sum_{1 \leq k \leq N, ~k \not\in A}  p(n_{\sigma(k)}x) \right| > \ve^{1/6} \right) \nonumber\\
& \leq & P_1^* ([-y-3\ve^{1/6},y+3\ve^{1/6}]) + H(\R \backslash [-\ve^{1/6},\ve^{1/6}]) + C_4 \ve^{2/3} + C_9 \ve^{2/3}\nonumber\\
& \leq & P_2^* ([-y-3\ve^{1/6},y+3\ve^{1/6}]) + 2 C_8 (y+3\ve^{1/6}) \ve^{1/2} \nonumber\\
& & ~+ 4 (y+3 \ve^{1/6}) \ve^{-1/3} e^{-\ve^{-1/3}/2} + \frac{2}{\sqrt{2 \pi \tau^2}} \int_{\ve^{1/6}}^\infty e^{-u^2/(2 \tau^2)} ~du + (C_4 + C_9) \ve^{2/3}\nonumber\\
& \leq & P_2 ([-y-4\ve^{1/6},y+4\ve^{1/6}]) + 2 C_8 (y+3\ve^{1/6}) \ve^{1/2} \nonumber\\
& & ~+ 4 (y+3\ve^{1/6}) \ve^{-1/3} e^{-\ve^{-1/3}/2} + \frac{4}{\sqrt{2 \pi \tau^2}} \int_{\ve^{1/6}}^\infty e^{-u^2/(2 \tau^2)} ~du  + (C_4 + C_9) \ve^{2/3}\nonumber\\
& \leq & P_2 ([-y,y]) \nonumber\\
& &  ~+ \frac{2}{\sqrt{2 \pi \gamma_N^2}} \int_{y}^{y+4\ve^{1/6}} e^{-u^2/(2 \gamma_N^2)} ~du  \label{gammaab}\\
& & ~ +2 C_8 (y+3\ve^{1/6}) \ve^{1/2} + 4 (y+3\ve^{1/6}) \ve^{-1/3} e^{-\ve^{-1/3}/2} \label{gammaab3}\\
& & ~+ \frac{4}{\sqrt{2 \pi \ve^{2/3}}} \int_{\ve^{1/6}}^\infty e^{-u^2/(2 \ve^{2/3})} ~du + (C_4 + C_9) \ve^{2/3}.\label{gammaab2}
\end{eqnarray}
It is clear that for $\gamma_N \to 0$ the limit distribution of $N^{-1/2} \sum_{k=1}^N f(n_{\sigma(k)}x)$ is the distribution  concentrated at 0. Now assume
$$
\gamma_N \to \gamma \qquad \textrm{as} \quad N \to \infty \qquad \textrm{for some} \quad \gamma>0.
$$
Then $\liminf_{N \to \infty} \gamma_N > 0$, and since $\ve$ can be chosen arbitrarily, the value of (\ref{gammaab}), (\ref{gammaab2}) and (\ref{gammaab3}) can also be made arbitrarily small. This means that for any given $\hat{\ve}>0$ and sufficiently large $N$
\begin{eqnarray}
& & \p \left( x:~N^{-1/2} \sum_{k=1}^N f(n_{\sigma(k)}x) \in [-y,y]\right) \label{newlab}\\
& \leq & \frac{1}{\sqrt{2 \pi \gamma_N^2}} \int_{-x}^x e^{-u^2/(2 \gamma_N^2)}~du + \hat{\ve} \nonumber\\
& \leq & \frac{1}{\sqrt{2 \pi \gamma^2}} \int_{-x}^x e^{-u^2/(2 \gamma^2)}~du + 2 \hat{\ve}. \nonumber
\end{eqnarray}
In the same way we can get a lower bound for (\ref{newlab}), proving the first part of Theorem \ref{th3}. On the other hand, the sequence $\gamma_N$ is bounded since $(n_k)_{k \geq 1}$ is lacunary and thus if $\gamma_N$ is not convergent, it has at least two accumulation points $\gamma^{(1)} \neq \gamma^{(2)}$.  Thus there exist two subsequences of $\N$ along which $N^{-1/2} \sum_{k=1}^N f(n_{\sigma(k)}x)$ converges to different distributions $\mathcal{N}(0,\gamma^{(1)}),~\mathcal{N}(0,\gamma^{(2)})$ (one of them may be the measure concentrated at 0, in case $\gamma^{(1)}=0$ or $\gamma^{(2)}=0$). But then  $N^{-1/2} \sum_{k=1}^N f(n_{\sigma(k)}x)$ does not have a single limit distribution as $N \to \infty$.\\

\subsection{Proof of Theorem \ref{th3bp}.} \label{sectth3bpproof}
In this section we prove Theorem \ref{th3bp} stating that the Diophantine condition $\mathbf{D}_2$ in Theorem \ref{th3} is best possible. In other words we show that if a lacunary sequence $(n_k)_{k \geq 1}$ does not satisfy condition $\mathbf{D}_2$, then there exist a trigonometric polynomial $f$ satisfying (\ref{fcond}) and a permutation $\sigma:~\N \to \N$ such that
$$
N^{-1/2} \sum_{k=1}^N f(n_{\sigma(k)}x)
$$
has a non-Gaussian limiting distribution.\\

Let $(n_k)_{k \geq 1}$ be given, and assume that condition $\mathbf{D}_2$ does not hold for this sequence. Then there exist integers $0<a\leq b$ and a sequence $(c_\ell)_{\ell \geq 1}$ of different  positive integers such that
$$
\# \left\{ k_1,k_2 \in \N:~ a n_{k_1} - b n_{k_2} = c_\ell\right\} \to \infty \qquad \textrm{as} \qquad \ell \to \infty.
$$
Throughout this section we will assume that $a<b$; the case $a=b$ can be handled in a similar way, with some minor changes.\\

We divide the set of positive integers into consecutive blocks $\Delta_1,\Delta_2,\dots,\Delta_m,\dots$ of lenghts $2^{2^1}, 2^{2^2}, \dots, 2^{2^m},\dots$; we will write $|\Delta_m|$ for the number of elements of a block $\Delta_m$. We can find integers $c_m,~m \geq 1$ and a permutation $\sigma:~\N \to \N$ such that the following holds:
\begin{eqnarray}
& & \sigma(k+1) > \sigma(k), \qquad\qquad k \geq 1 \nonumber\\
& & \sigma(2k) > 2b \, \sigma(2k-1), \qquad\qquad k \geq 1 \nonumber\\
& & a n_{\sigma(2k)} - b n_{\sigma(2k-1)} = c_m, \qquad 2k,2k-1 \in \Delta_m, ~m \geq 1. \label{diocm}
\end{eqnarray}
Let
$$
f(x)=\cos 2 \pi a x + \cos 2 \pi b x.
$$
To show that
$$
N^{-1/2}\sum_{k=1}^N f(n_{\sigma(k)} x)
$$
does not converge to a Gaussian distribution, we will show that
$$
\frac{\sum_{k \in \Delta_m} f(n_{\sigma(k)}x)}{\sqrt{|\Delta_m|}}
$$
has a non-Gaussian limit distribution as $m \to \infty$. Since the lengths of the block $\Delta_m$ dominates $|\Delta_1| + \dots + |\Delta_{m-1}|$, this means that
$$
\frac{\sum_{m=1}^M \sum_{k \in \Delta_m} f(n_{\sigma(k)}x)}{\sqrt{\sum_{m=1}^M |\Delta_m|}}
$$
also has a non-Gaussian limiting distribution. \\

\begin{lmm} \label{lemmanon}
$$
\left| \E \left( \exp \left( \frac{i s \sum_{k \in \Delta_m} f(n_{\sigma(k)} x)}{\sqrt{|\Delta_m|}} \right)\right) - \E \left(e^{-s^2 (1+\cos 2 \pi x)/2}\right)\right| \ll |\Delta_m|^{-1/8},
$$
for all $s \in [-(\log |\Delta_m|)^{1/8},-(\log |\Delta_m|)^{1/8}]$.
\end{lmm}

\noindent
\emph{Proof.~} We write $(\nu_k)_{k \geq 1}$ for $(n_{\sigma(k)})_{k \geq 1}$ and define
$$
\eta_m = \frac{\sum_{k \in \Delta_m} f(\nu_k x)}{\sqrt{|\Delta_m|}}.
$$
For different $k_1,k_2$, for which $2k_1-1,2k_2-1 \in\Delta_m$, the functions
$$
1 + \frac{i s \left(f(\nu_{2 k_1-1}x)+f(\nu_{2k_1}x) \right)}{\sqrt{2}} \qquad \textrm{and} \qquad 1 + \frac{i s \left(f(\nu_{2 k_2-1}x)+f(\nu_{2k_2}x)\right)}{\sqrt{2}}
$$
are orthogonal. This means, writing
$$
\alpha_m(s)=\prod_{k\geq 1:~2k-1 \in \Delta_m} \left(1 + \frac{i s \left(f(\nu_{2 k-1}x)+f(\nu_{2k}x) \right)}{\sqrt{2 |\Delta_m|}} \right),
$$
we have
$$
\E \Big(\alpha_m(s)\Big)=1.
$$
Using again
(\ref{win}), we get
\begin{eqnarray*}
e^{is \eta_m} & = & \prod_{k\geq 1:~2k-1 \in \Delta_m} \exp\left(\frac{i s \left(f(\nu_{2 k-1}x)+f(\nu_{2k}x) \right)}{\sqrt{2|\Delta_m|}} \right) \\
& = & \alpha_m(s) ~\exp \left(\sum_{k\geq 1:~2k-1 \in \Delta_m} \frac{-s^2 \left(f(\nu_{2 k-1}x)+f(\nu_{2k}x) \right)^2}{2|\Delta_m|} \right) ~ W_m(s),
\end{eqnarray*}
where
\begin{eqnarray*}
W_m(s)= \exp \left(\sum_{k\geq 1:~2k-1 \in \Delta_m} w \left(\frac{i s \left(f(\nu_{2 k-1}x)+f(\nu_{2k}x) \right)}{\sqrt{2|\Delta_m|}}\right)\right).
\end{eqnarray*}
By (\ref{diocm}) we have
\begin{eqnarray*}
& & \sum_{k\geq 1:~2k-1 \in \Delta_m} \left(f(\nu_{2 k-1}x)+f(\nu_{2k}x) \right)^2 \\
& = & \sum_{k\geq 1:~2k-1 \in \Delta_m} \left( \cos 2 \pi a \nu_{2k-1} x + \cos 2 \pi b \nu_{2k-1}x + \cos 2 \pi a \nu_{2k} x + \cos 2 \pi b \nu_{2k}x\right)^2 \\
& = & |\Delta_m| + \left(\sum_{k\geq 1:~2k-1 \in \Delta_m} \frac{\cos 2 \pi (a \nu_{2k} - b \nu_{2k_1}) x}{2} \right) + R_m(x) \\
& = & |\Delta_m| +  \left(\sum_{k\geq 1:~2k-1 \in \Delta_m} \frac{\cos 2 \pi c_m x}{2} \right) + R_m(x) \\
& = & |\Delta_m| \left( 1 + \frac{\cos 2 \pi c_m x}{4} \right) + R_m(x),
\end{eqnarray*}
where $R_m(x)$ is a sum of cosine functions with coefficients 1/2 and frequencies
\begin{eqnarray*}
& & a \nu_{2k-1}, ~b \nu_{2k-1},~a \nu_{2k},~b \nu_{2k}, ~a \nu_{2k-1} \pm b \nu_{2k-1},~a \nu_{2k-1} \pm a \nu_{2k}, ~a \nu_{2k-1} \pm b \nu_{2k},\\
& & b \nu_{2k-1} + a \nu_{2k},~b \nu_{2k-1} \pm b \nu_{2k},~a \nu_{2k} \pm b \nu_{2k},
\end{eqnarray*}
where $k$ runs through the set $\{k \geq 1: ~2k-1 \in \Delta_m\}$. If we write
\begin{eqnarray*}
R_m^{(1)} (x) & = & \frac{1}{2} \sum_{k \geq 1: ~2k-1 \in \Delta_m} \cos 2 \pi a \nu_{2k-1} x, \\
R_m^{(2)} (x) & = & \frac{1}{2} \sum_{k \geq 1: ~2k-1 \in \Delta_m} \cos 2 \pi b \nu_{2k-1} x, \\
\vdots \\
R_m^{(14)} (x) & = & \frac{1}{2} \sum_{k \geq 1: ~2k-1 \in \Delta_m} \cos 2 \pi (a \nu_{2k} + b \nu_{2k})x, \\
R_m^{(15)} (x) & = & \frac{1}{2} \sum_{k \geq 1: ~2k-1 \in \Delta_m} \cos 2 \pi (a \nu_{2k} - b \nu_{2k})x,
\end{eqnarray*}
the function $R_m(x)$ can be split into 15 lacunary cosine series, each consisting of $\Delta_m/2$ elements. Observing
\begin{eqnarray} \label{alphaest}
|\alpha_m(s)| \leq \prod_{k \geq 1:~2k-1 \in \Delta_m} \left(1+\frac{4s^2}{2|\Delta_m|}\right)^{1/2} \leq e^{s^2/2}
\end{eqnarray}
and, in view of (\ref{win}),
\begin{eqnarray}
|W_m(s)-1| & \leq & \left|\exp \left( \frac{|\Delta_m|}{2} ~ \left| \frac{2s}{\sqrt{2 |\Delta_m|}} \right|^3  \right) -1\right| \nonumber\\
& \leq & e^{4 |\Delta_m|^{-1/2}} - 1\nonumber\\
& \leq & \frac{8}{|\Delta_m|^{1/2}}. \nonumber\label{ewest}
\end{eqnarray}
We have
\begin{eqnarray}
& & \left| \E e^{i s \eta_m} - \E \exp \left(-\frac{s^2}{2} \left( 1+\frac{\cos 2 \pi x}{4} \right)  \right) \right| \label{wehav}\\
& = & \left| \E \left( \alpha_m(s) ~\exp \left(-\frac{s^2}{2} \left(1+\frac{\cos 2 \pi c_m x}{4}\right) \right)~\exp \left(\frac{-s^2R_m(x)}{2 |\Delta_m|} \right) ~W_m(s) \right) \right. \nonumber\\
& & \quad \left. - \E \exp \left(-\frac{s^2}{2} \left(1+\frac{\cos 2 \pi c_m x}{4} \right) \right) \right| \nonumber\\
& \leq & \left| \E \left( \alpha_m(s) ~\exp \left(-\frac{s^2}{2} \left(1+\frac{\cos 2 \pi c_m x}{4}\right) \right) \times \right.\right.\label{integ1}\\
& & \quad ~\times \left.\left.\left( \exp \left(\frac{-s^2R_m(x)}{2 |\Delta_m|} \right) ~W_m(s) -1 \right) \right) \right| \label{integ2}\\
& & + \left| \E \left( \left(\alpha_m(s)-1\right) ~\exp \left(-\frac{s^2}{2} \left(1+\frac{\cos 2 \pi c_m x}{4} \right) \right) \right)  \right| \label{integ4}
\end{eqnarray}
In view of (\ref{alphaest}), the term in lines (\ref{integ1}), (\ref{integ2}) is at most
\begin{eqnarray}
& & \max_{x \in (0,1)} \left(\exp\left(-\frac{s^2 \cos 2 \pi c_m x}{8} \right)\right) ~\E \left|\exp \left(\frac{-s^2R_m(x)}{2 |\Delta_m|} \right) ~W_m(s) -1 \right| \nonumber\\
& \leq & e^{s^2/8} ~\E \left|\exp \left(\frac{-s^2R_m(x)}{2 |\Delta_m|} \right) ~W_m(s) -1 \right|. \label{wehav2}
\end{eqnarray}
It is easy to see that
$$
\|R_m\| \leq 15 |\Delta_m|
$$
and therefore
$$
\p \left\{ |R_m| \geq 15 |\Delta_m|^{2/3} \right\} \leq \Delta_m^{-1/3}.
$$
Thus
\begin{eqnarray}
& & \E \left|\exp \left(\frac{-s^2R_m(x)}{2 |\Delta_m|} \right) ~W_m(s) -1 \right| \label{wehav3}\\
& \leq & \left(1+\frac{15 s^2}{|\Delta_m|^{1/3}} \right) \left(1+\frac{8}{|\Delta_m|^{1/2}} \right) -1 + \frac{15s^2 |\Delta_m|}{4 |\Delta_m|^{4/3}} \nonumber\\
& \ll & s^2 |\Delta_m|^{-1/3},\label{wehav4}
\end{eqnarray}
where the implied constant does not depend on $m,s$.\\

The function $\alpha_m(s)$ is a sum of the constant term $1$ plus at most $2^{2|\Delta_m|}$ cosine functions with coefficients at most $1$ (provided $|s| \leq |\Delta_m|^{1/2}$) and frequencies at least $2^{4 |\Delta_m|}$. Thus (\ref{integ4}) is at most
\begin{eqnarray}
& & 2^{2 |\Delta_m|} 2^{-4|\Delta_m|} \max_{x \in (0,1)} \left| \frac{d}{dx}~ \exp \left(-\frac{s^2}{2} \left(1+\frac{\cos 2 \pi c_m x}{4} \right) \right) \right| \nonumber\\
& \leq & 2^{-2 |\Delta_m|} ~e^{5s^2/8} ~\frac{\pi s^2 c_m}{4} \label{wehav5}.
\end{eqnarray}
Combining (\ref{wehav2}), (\ref{wehav3}), (\ref{wehav4}) and (\ref{wehav5}) we see that (\ref{wehav}) is at most
$$
\ll e^{s^2/8} s^2 |\Delta_m|^{-1/3} + 2^{-2 |\Delta_m|} ~e^{5s^2/8} ~\frac{\pi s^2 c_m}{4}.
$$
In particular, since we assumed $|s| \leq (\log |\Delta_m|)^{1/8}$, we get that the expression (\ref{wehav}) is at most
$$
\ll |\Delta_m|^{-1/8},
$$
proving Lemma \ref{lemmanon}.\\

\noindent
\emph{Proof of Theorem \ref{th3bp}.}~The proof of Theorem \ref{th3bp} can be obtained from Lemma \ref{lemmanon} like the proof of Theorem \ref{th3} was obtained from Lemma \ref{lemmaclt} in Section \ref{sectth3proof}. The normed partial sums
$$
|\Delta_m|^{-1/2} \sum_{k \in \Delta_m} f(n_{\sigma(k)} x)
$$
have a limiting distribution which is non-Gaussian, and since the set $\Delta_m$ is much larger than the union of the sets $\Delta_1, \dots, \Delta_{m-1}$, the same distribution is the limit distribution of
$$
N^{-1/2} \sum_{k=1}^N f(n_{\sigma(k)}x).
$$
This proves Theorem \ref{th3bp}.


\end{document}